\documentclass{article}[11pt]
\usepackage{hyperref}
\usepackage{amssymb}
\usepackage{amsfonts}
\usepackage{amsmath}
\usepackage{amsthm}
\usepackage{mathtools}
\usepackage{array}
\usepackage{amsmath}
\usepackage{amsfonts}
\usepackage{amssymb}
\usepackage{pgf}
\usepackage{tikz}
\usepackage{tkz-tab}
\usetikzlibrary{shapes,snakes,arrows,backgrounds}
\usetikzlibrary{scopes,svg.path,shapes.geometric,shadows}
\usepackage{graphicx}
\usepackage{xcolor}
\usepackage{tabu}

\begin{document}
\title{A Simple Proof of the Riemann's Hypothesis}
\author{Charaf ECH-CHATBI  \thanks{Email: charaf.chatbi@gmail.com. 
The opinions of this article are those of the author and do not reflect in any way the views or business of his employer.} }
\date{23 October 2023} 
\maketitle

\DeclarePairedDelimiter\abs{\lvert}{\rvert}%
\DeclarePairedDelimiter\norm{\lVert}{\rVert}%

\newcommand{\dd}{\mathrm{d}}
\newcommand{\ee}{\mathrm{e}}
\newcommand{\bbbp}{\mathbb{P}}
\newcommand{\bbbe}{\mathbb{E}}
\newcommand{\bbbr}{\mathbb{R}}
\newcommand{\QQ}{\mathbb{Q}}
\newcommand{\EE}{\mathbb{E}}
\newcommand{\PP}{\mathbb{P}}
\newcommand{\SRR}{$\sigma_R $ }
\newcommand{\ERR}{$\overline{R} $}
\newcommand{\CRO}{\textsc{CRO} }
\newcommand{\I}{\mathrm{1\!\!1}}

\newcommand*{\QEDA}{\null\nobreak\hfill\ensuremath{\blacksquare}}%
\newcommand*{\QEDB}{\null\nobreak\hfill\ensuremath{\square}}%

\newtheorem{theorem}{Theorem}[section]
\newtheorem{corollary}{Corollary}[theorem]
\newtheorem{lemma}[theorem]{Lemma}
\newtheorem*{remark}{Remark}
\newtheorem*{strategy}{Proof Strategy}
\newtheorem*{conclusion}{Conclusion}

%

\makeatletter
\newcommand{\subjclass}[2][1991]{%
  \let\@oldtitle\@title%
  \gdef\@title{\@oldtitle\footnotetext{#1 \emph{Mathematics subject classification.} #2}}%
}
\newcommand{\keywords}[1]{%
  \let\@@oldtitle\@title%
  \gdef\@title{\@@oldtitle\footnotetext{\emph{Key words and phrases.} #1.}}%
}
\makeatother

 \abstract{ We present a simple proof of the Riemann's Hypothesis (RH) where only undergraduate mathematics is needed.} \\


    

\section{The Riemann Hypothesis}

\subsection{The importance of the Riemann Hypothesis}

The prime number theorem gives us the average distribution of the primes. The Riemann hypothesis tells us about the deviation from the average. Formulated in Riemann's 1859 paper[1], it asserts that all the 'non-trivial' zeros of the zeta function are complex numbers with real part 1/2.

\subsection{Riemann Zeta Function}

For a complex number $s$ where $\Re(s) > 1$, the Zeta function is defined as the sum of the following series:  \begin{eqnarray}   \zeta(s) = \sum^{+\infty}_{n = 1} \frac{1}{n^s} \end{eqnarray}
In his 1859 paper[1], Riemann went further and extended the zeta function $\zeta(s)$, by analytical continuation, to an absolutely convergent function in the half plane $\Re(s) > 0$, minus a simple pole at s = 1:
\begin{eqnarray}   \zeta(s) =  \frac{s}{s - 1} - s \int^{+\infty}_{1} \frac{ \{x\} }{ x^{s + 1} } dx  \end{eqnarray}
Where $ \{x\} = x -  [x] $ is the fractional part and $[x]$ is the integer part of $x$. 
Riemann also obtained the analytic continuation of the zeta function to the whole complex plane.

Riemann[1] has shown that Zeta has a functional equation\footnote{This is slightly different from the functional equation presented in Riemann's paper[1]. This is a variation that is found everywhere in the litterature[2,3,4]. Another variant using the $\cos$: \begin{eqnarray} \zeta(1 - s)   =  2^{1 - s} \pi^{-s} \cos \big(  \frac{\pi s}{2} \big) \Gamma(s) \zeta(s) \end{eqnarray} }
\begin{eqnarray} \zeta(s)   =  2^{s} \pi^{s-1} \sin \big(  \frac{\pi s}{2} \big) \Gamma(1-s) \zeta(1-s) \end{eqnarray} 
Where $\Gamma(s)$ is the Gamma function. Using the above functional equation, Riemann has shown that the non-trivial zeros of $\zeta$ are located symmetrically with respect to the line $\Re(s) = 1/2$, inside the critical strip $0 < \Re(s) < 1$.  Riemann has conjectured that all the non trivial-zeros are located on the critical line $\Re(s) = 1/2$. 
In 1921, Hardy $\&$ Littlewood[2,3, 6] showed that there are infinitely many zeros on the critical line. In 1896, Hadamard and De la Vall\'ee Poussin[2,3] independently proved that $\zeta(s)$ has no zeros of the form $s = 1 + \it{i}t$ for $t  \in \mathbb{R}$. Some of the known results[2, 3] of $\zeta(s)$ are as follows:
\begin{itemize}

\item $\zeta(s)$ has no zero for $\Re(s) > 1$.
\item  $\zeta(s)$ has no zero of the form $s = 1 + \it{i}\tau$. i.e. $\zeta(1+ i \tau) \neq 0$, $\forall \, \tau$.
\item $\zeta(s)$ has a simple pole at $s = 1$ with residue 1.
\item $\zeta(s)$ has all the trivial zeros at the negative even integers $s = -2k$, $k \in \mathbb{N^*}$.
\item The non-trivial zeros are inside the critical strip: i.e. $0 < \Re(s) < 1$.
\item If $\zeta(s) = 0$, then $1-s$, $\bar{s}$ and $1 - \bar{s}$ are also zeros of $\zeta$: i.e. $ \zeta(s) =  \zeta(1 - s) =  \zeta( \bar{s} ) =  \zeta( 1 - \bar{s} )  = 0$.

\end{itemize}

Therefore, to prove the “Riemann Hypothesis” (RH), it is sufficient to prove that $\zeta$ has no zero on the right hand side $1/2 < \Re(s) < 1$ of the critical strip.

\subsection{Proof of the Riemann Hypothesis}

Let's take a complex number $s$ such that $s = \sigma + \it{i} \tau $. Unless we explicitly mention otherwise, let's suppose that $ 0 < \sigma < 1$, $\tau > 0$  and $\zeta (s) = 0$.  \\

We have from the Riemann's integral above: \begin{eqnarray}  \zeta(s) =  \frac{s}{s - 1} - s \int^{+\infty}_{1} \frac{ \{x\} }{ x^{s + 1} } dx  \end{eqnarray}
We have $s \neq 1$, $s \neq 0$ and $\zeta(s) = 0$, therefore:  \begin{eqnarray}  \frac{1}{s - 1} = \int^{+\infty}_{1} \frac{ \{x\} }{ x^{s + 1} } dx  \end{eqnarray}
Or in other terms: \begin{eqnarray}  \frac{1}{\sigma + \it{i} \tau - 1} = \int^{+\infty}_{1} \frac{ \{x\} }{ x^{\sigma + \it{i} \tau + 1} } dx  \end{eqnarray}

Let's denote the following functions:
\begin{eqnarray}  
\epsilon(x) & = &  \{x\}   \\
\phi(x) & = &  \{ x\}  \big( 1 - \{ x \} \big) \\
\Psi(x) & = & \int^{x}_{1} \!\!\! dx \, \epsilon(x) \\
K(\sigma, \tau) & = &   \frac{\tau}{(1-\sigma)^2 + \tau^2}   \hspace{1cm} 
\end{eqnarray}

To continue, we will prove the following lemmas.

\begin{lemma}
Let's consider two variables $\sigma$ and $\tau$ such that $\sigma >0$ and $\tau >0$. Let's define two integrals $I(a, \sigma, \tau)$ and $J(a, \sigma, \tau)$ as follows: 
\begin{eqnarray} I(a, \sigma, \tau) & = &  \int^{a}_{1}  \frac{  \sin (\tau \ln(x)) }{ x^{\sigma} }   dx \hspace{1.5cm}   \\
J(a, \sigma, \tau) & = &  \int^{a}_{1}  \frac{  \cos (\tau \ln(x)) }{ x^{\sigma} }   dx   \hspace{1.5cm}  
\end{eqnarray}  
Therefore
\begin{eqnarray} I(a, \sigma, \tau) & = &  K(\sigma, \tau) \Big(   1 - \frac{  \cos (\tau \ln(a))   }{ a^{\sigma - 1}  }    -  \frac{ ( \sigma - 1 )}{ \tau  }    \frac{  \sin (\tau \ln(a))    }{ a^{\sigma - 1}    }       \Big)   \hspace{0.6cm} \\
J(a, \sigma, \tau) & = &  K(\sigma, \tau) \Big(  \frac{ (\sigma - 1)}{\tau} -  \frac{(\sigma - 1)}{\tau} \frac{  \cos (\tau \ln(a))   }{ a^{\sigma - 1}  }    +     \frac{  \sin (\tau \ln(a))    }{ a^{\sigma - 1}    }   \Big)   \hspace{0.6cm} 
\end{eqnarray}

\end{lemma}

\begin{proof}

Let's consider two variables $\sigma$ and $\tau$ such that $\sigma >0$ and $\tau >0$. Let's take $a > 1$.
\begin{eqnarray} I(a, \sigma, \tau) & = &  \int^{a}_{1}  \frac{  \sin (\tau \ln(x)) }{ x^{\sigma} }   dx  \\
& = & \int^{\ln(a)}_{0}   \sin (\tau x)   e^{(1- \sigma) x} dx \\
& = & K(\sigma, \tau) \Big(   1 - \frac{  \cos (\tau \ln(a))   }{ a^{\sigma - 1}  }    -  \frac{ ( \sigma - 1 )}{ \tau  }    \frac{  \sin (\tau \ln(a))    }{ a^{\sigma - 1}    }       \Big)   \hspace{1.2cm} \end{eqnarray}

And the same for $J(a, \sigma, \tau)$ for $a >0$:
\begin{eqnarray} J(a, \sigma, \tau) & = &  \int^{a}_{1}  \frac{  \cos (\tau \ln(x)) }{ x^{\sigma} }   dx  \\
& = & \int^{\ln(a)}_{0}   \cos (\tau x)   e^{(1- \sigma) x} dx \\
& = & K(\sigma, \tau) \Big(  \frac{ (\sigma - 1)}{\tau} -  \frac{(\sigma - 1)}{\tau} \frac{  \cos (\tau \ln(a))   }{ a^{\sigma - 1}  }    +     \frac{  \sin (\tau \ln(a))    }{ a^{\sigma - 1}    }        \Big)   \hspace{0.6cm} \end{eqnarray}

\end{proof}

\begin{lemma}  
The function $\epsilon(x)$ is piecewise continuous on $[0, +\infty)$ and its primitive function $\Psi(x)$ is defined as follows:
\begin{eqnarray}  \Psi(x) & = &  \frac{1}{2} \Big( x - 1 - \phi(x) \Big)   \end{eqnarray}
Let's consider two variables $\sigma$ and $\tau$ such that $0 < \sigma < 1$ and $\tau >0$ such that $s = \sigma + \textit{i} \tau$ is a zeta zero. 
Therefore:  \begin{eqnarray} 
\int^{+\infty}_{1}   \hspace{-0.25cm} dx \frac{  \Psi(x)  }{  \, x^{2 + s}  }     & = &  \frac{1}{  (s - 1) (1 + s)  } \hspace{1cm} 
\end{eqnarray}  
And \begin{eqnarray} 
 \int^{+\infty}_{1}   \hspace{-0.25cm} dx \frac{ \phi(x)  }{  \, x^{2 + s }  }   & = &  \frac{1  }{  s (1 - s)   } \hspace{1cm} 
\end{eqnarray}  
\end{lemma}

\begin{proof}
We will need the function $\phi$ as you will see later that we need a continuous function instead of a piecewise one like the function $\epsilon$. \\
Let's take $x > 1$ a real number. Let's denote $n_x = \lfloor{x}\rfloor$ be the integer part of $x$. We have $n_x = x - \{ x \}$. Therefore, we can write the following:
 \begin{eqnarray}  \Psi(x) & = & \int^{x}_{1} \epsilon(t) dt  \hspace{1cm} \\
& = & \sum^{n_x-1}_{n = 1}  \int^{n+1}_{n} \{ t \} dt +   \int^{x}_{n_x} \{ t \} dt \hspace{1cm} \\ 
& = & \sum^{n_x-1}_{n = 1}  \int^{n+1}_{n} (t - n) dt+   \int^{x}_{n_x} (t - n_x) dt \hspace{1cm} \\ 
& = & \sum^{n_x-1}_{n = 1}  \frac{1}{2}+   \frac{1}{2} (x-n_x)^2 \hspace{1cm} \\
& = & \frac{1}{2} \Big( n_x - 1  +    \{ x \}^2 \Big) \hspace{1cm} \\
& = & \frac{1}{2} \Big( x - 1  - \{ x \}+    \{ x \}^2 \Big) \hspace{1cm} 
  \end{eqnarray}
This prove the equation (23).
\QEDA

Let's prove the second point of the lemma. Let's define the integral $I_{\epsilon}(s)$ as follows: 
\begin{eqnarray}   
I_{\epsilon}(s) & = &   \int^{+\infty}_{1}   \hspace{-0.25cm} dx \frac{ \epsilon(x)  }{  \, x^{1 + s }  }   
\end{eqnarray}  

The function $x \to \frac{  \epsilon(x)  }{  \, x^{1 + s }  } $ is integrable on $[1, +\infty)$ and thanks to the integration by parts, we can write the following:
\begin{eqnarray}   I_{\epsilon}(s) & = &   \int^{+\infty}_{1}   \hspace{-0.25cm} dx \frac{ \epsilon(x)  }{  \, x^{1 + s}  }   \hspace{1cm}  \\
& = & \Big[ \frac{  \Psi(x)  }{  \, x^{1 + s }  }   \Big]^{+\infty}_{x=1}  + (1 + s)  \int^{+\infty}_{1}   \hspace{-0.25cm} dx \frac{  \Psi(x)  }{  \, x^{2 + s }  }   \hspace{1cm} \\
& = & (1 + s) \int^{+\infty}_{1}   \hspace{-0.25cm} dx \frac{   \Psi(x) }{  \, x^{2 + s }  } 
 \end{eqnarray}  

Since, $s$ is a zeta zero, from the equation (6), we have:

\begin{eqnarray}   
I_{\epsilon}(s) & = &    \frac{1}{s-1}
\end{eqnarray}

Therefore

\begin{eqnarray}   
\int^{+\infty}_{1}   \hspace{-0.25cm} dx \frac{   \Psi(x) }{  \, x^{2 + s }  }  & = & \frac{1}{(s-1) (s+1)}
\end{eqnarray}

Thanks to equation (23), we can write:
\begin{eqnarray}  \hspace{-0cm} 
 \int^{+\infty}_{1}   \hspace{-0.25cm} dx \frac{   \phi(x)   }{  \, x^{2 + s }  }   & = &   \int^{+\infty}_{1}   \hspace{-0.25cm} dx \frac{ ( x - 1  )   }{  \, x^{2 + s }  }   - 2 \int^{+\infty}_{1}   \hspace{-0.25cm} dx \frac{   \Psi(x) }{  \, x^{2 + s }  }  \hspace{1cm} \\
 & = &  \frac{1}{ s  } - \frac{1}{s+1}   -    \frac{2}{ (s-1) (s+1)  }  \hspace{1cm} \\
& = &   \frac{1}{ s(1-s)  }   \hspace{1cm} 
 \end{eqnarray}  

 \QEDA
\end{proof}



\begin{lemma}
Let's consider two variables $\sigma$ and $\tau$ such that $0 < \sigma < 1$, $\tau >0$ and $s = \sigma + \textit{i} \tau$ is a zeta zero. 
Let's define the sequence of functions $\phi_n$ and $\psi_n$ over $[0,+\infty)$ and $\overline{\phi}_n$ over $[1, +\infty)$ such that $\phi_0(x) = \overline{\phi}_0(x) =  \phi(x) $ and for each $n \geq 1$, $\phi_n(0) = 0 $, and: 
\begin{eqnarray}  
\phi_{n+1}(x) & = & \frac{1}{x} \int^{x}_{0} dt \, \phi_{n}(t)  \hspace{1cm}  \text{for  $x >0$} \hspace{1cm}   \\
\overline{\phi}_{n+1}(x) & = & \frac{1}{x} \int^{x}_{1} dt \, \overline{\phi}_{n}(t)   \hspace{1cm}  \text{for  $x \geq 1$} \hspace{1cm}  \\
\psi_{n}(x) & = & \frac{x}{2^n} -  \frac{x^2}{3^n}   \hspace{1.7cm}  \text{for  $x \geq 0$} \hspace{1cm}   
\end{eqnarray} 

Therefore:
\begin{enumerate}
\item For each $n \geq 1$:

\begin{eqnarray}  
 \int^{+\infty}_{1}   \hspace{-0.25cm} dx \frac{ \phi_{n}(x)    }{  \, x^{2 + s }  }  & = &   \frac{1}{s \, 2^{n}} + \frac{1}{(1-s) \, 3^{n}}       \hspace{1cm}  
\end{eqnarray}

\item For each $n$:
\begin{eqnarray}  
 \int^{+\infty}_{1}   \hspace{-0.25cm} dx \frac{ \phi_{n}(x)    }{  \, x^{3 - s }  }  & = &   \frac{1}{(1-s) \, 2^{n}} + \frac{1}{s \, 3^{n}}       \hspace{1cm}  
\end{eqnarray}

\item For each $x \geq 1$:
\begin{eqnarray}  \phi_{n}(x) =  \overline{\phi}_{n}(x) + \frac{1}{ 2^n \, x } \sum^{n-1}_{k = 0} \frac{ 2^k  \ln^k(x)  }{k!} - \frac{1}{ 3^n \, x } \sum^{n-1}_{k = 0} \frac{ 3^k  \ln^k(x)  }{k!}  \end{eqnarray} 

\item For each $n$ and $x \geq 1$:
\begin{eqnarray}  
  \overline{\phi_{n}}(x) & = &  \frac{1}{(n-1)!}  \frac{1}{x}    \int^{x}_{1} dt\,  \overline{\phi_{0}}(t)  \big( \ln(\frac{x}{t}) \big)^{n-1}   \hspace{1cm}
\end{eqnarray}

\item For each $x \geq 0$:
\begin{eqnarray}  \lim_{n \to +\infty } 2^n  \,  \phi_{n}(x) & = & x \end{eqnarray} 

\item For each $x \geq 0$:
\begin{eqnarray}  0  \leq  \,\,  2^n  \phi_{n}(x) & \leq & x \end{eqnarray} 

\item For each $x \geq 0 $:
\begin{eqnarray} 
\psi_n(x) & = &  \overline{\psi_{n}}(x)   + \frac{1}{ 2^n \, x } \sum^{n-1}_{k = 0} \frac{ 2^k  \ln^k(x)  }{k!} - \frac{1}{ 3^n \, x } \sum^{n-1}_{k = 0} \frac{ 3^k  \ln^k(x)  }{k!}    \hspace{1cm}
\end{eqnarray} Where 
\begin{eqnarray} \overline{\psi_{n}}(x)  & = & \frac{ 1 }{(n-1)!}  \frac{1}{x}  \int^{x}_{1} dt\,  (t - t^2)  \big( \ln(\frac{x}{t}) \big)^{n-1} \end{eqnarray}

\end{enumerate}
\end{lemma}

\begin{proof}
We have $s= \sigma + \textit{i} \tau$ a $\it{zeta}$ zero. The lemma 1.2 calculates the integral $A(s)$ as follows:
\begin{eqnarray} A(s) & = & \int^{+\infty}_{1}   \hspace{-0.25cm} dx \frac{  \phi_0(x)   }{  \, x^{2 + s }  }  = \frac{1}{s(1-s)}  \hspace{1cm}     \end{eqnarray} 

We use the integration by parts to write the following:
\begin{eqnarray} \hspace{-1cm} 
 A(s) & = &  \int^{+\infty}_{1}   \hspace{-0.25cm} dx \frac{ \phi_0(x)    }{  \, x^{2 + s }  }    \hspace{0.7cm}  \\
& = & \Bigg[  \frac{ 1 }{  \, x^{2 + s }  }  \int^{x}_{0} \!\! dx \, \phi_0(x)    \Bigg]^{+\infty}_{1} + (2 + s)  \int^{+\infty}_{1}   \hspace{-0.25cm} dx \frac{ \phi_1(x)    }{  \, x^{2 + s }  }    \hspace{0.7cm}  \\
& = & - \int^{1}_{0} \!\! dx \, \phi_0(x)   + (2 + s)  \int^{+\infty}_{1}   \hspace{-0.25cm} dx \frac{ \phi_1(x)    }{  \, x^{2 + s }  }    \hspace{0.7cm}  \\
& = & - \int^{1}_{0} \!\! dx \, \phi_0(x)  - (2 + s) \int^{1}_{0} \!\! dx \, \phi_1(x)   + (2 + s)^2  \int^{+\infty}_{1}   \hspace{-0.25cm} dx \frac{ \phi_2(x)    }{  \, x^{2 + s }  }    \hspace{0.7cm}  \\
& ... & \\
& = &   - \sum^{n}_{k=0} (2 + s)^k \int^{1}_{0}   \hspace{-0.25cm} dx \, \phi_{k}(x)  + (2 + s)^{n+1} \int^{+\infty}_{1}   \hspace{-0.25cm} dx \frac{ \phi_{n+1}(x)    }{  \, x^{2 + s }  }    \hspace{0.7cm}  
\end{eqnarray}

We could do the above integration by parts because the functions $x \to \frac{\phi_k(x)}{x^{2+s}}$ are piecewise continuous\footnote{All the functions $x \to \frac{\phi_k(x)}{x^{2+s}}$ are continuous except when $k=0$ as the function $\phi_0$ is piecewise continuous.} and integrable over $[1,+\infty)$ as they are dominated by the function $x \to \frac{1}{x^{2+\sigma}}$ that is integrable over $[1,+\infty)$. In fact, we can prove that the functions $\phi_k$ are non-negative and bounded by $1$ as we can prove by recurrence that for each $k$ for each $x > 0$ that: \begin{eqnarray}  0 < \phi_{k+1}(x) = \frac{1}{x} \int^{x}_{0} dt \, \phi_{k}(t)   \leq    \frac{1}{x} \int^{x}_{0} 1 \, dt  = 1 \end{eqnarray}
This can be proven by recurrence using the fact that it is true for the initial case as we have for each $x \geq 0$:  \begin{eqnarray}  0 \leq \phi_0(x) =  \{ x\} (1 - \{ x\} )   &  \leq  &  1 \end{eqnarray}

Now, we need to calculate the integrals $I_k$ for $k \geq 0$: 
\begin{eqnarray} I_k & = &  \int^{1}_{0} \!\! dx \, \phi_k(x) \end{eqnarray}
For $x \in (0,1)$, we have: \begin{eqnarray}  \phi_0(x) =  x - x^2   \end{eqnarray}
And  \begin{eqnarray}  \phi_1(x) =  \frac{x}{2} - \frac{x^2}{3}   \end{eqnarray}

Therefore, we can write for each $k$ for $x \in (0,1)$: \begin{eqnarray}  \phi_k(x) =  \frac{x}{2^k} -  \frac{x^2}{ 3^k}   \end{eqnarray}
Therefore, for each $k \geq 1$: \begin{eqnarray}  I_k =  \frac{1}{2^{k+1}} - \frac{1}{ 3^{k+1}}   \end{eqnarray}

Therefore, we can conclude:

\begin{eqnarray}  \hspace{-1cm} 
 A(s) & = & \frac{1}{s(1-s)} +  (2 + s)^{n+1} \Bigg[ \int^{+\infty}_{1}   \hspace{-0.25cm} dx \frac{ \phi_{n+1}(x)    }{  \, x^{2 + s }  } - \Bigg( \frac{1}{s \, 2^{n+1}} + \frac{1}{(1-s) \, 3^{n+1}}   \Bigg) \Bigg]    \hspace{0.7cm} 
\end{eqnarray}

Since $2 + s \neq 0$ and $3 -s \neq 0$, we conculde the result of our lemma. \QEDA

The point 3) can be proved by recurrence. For $n = 0$. We have $ \phi_{0}(x) = \overline{\phi_0}(x)$. Let's assume that it is true till $n$ and let's prove for $n+1$. 
We have:
\begin{eqnarray}  
\phi_{n+1}(x)   & = &  \frac{1}{x} \int^{1}_{0}   \hspace{-0.25cm} dt \, \phi_n(t)   + \frac{1}{x} \int^{x}_{1}   \hspace{-0.25cm} dt \, \phi_n(t)    \hspace{1cm}  \\
& = &  \frac{ 1 }{ x } \big(\frac{ 1 }{ 2^{n +1 } } - \frac{ 1 }{ 3^{n +1 } } \big)  + \frac{1}{x \, 2^n}  \int^{x}_{1}   \hspace{-0.15cm} \, \frac{dt}{  t } \sum^{n-1}_{k = 0} \frac{ 2^k  \ln^k(t)  }{k!}   \hspace{0.5cm} \\
&  & - \frac{1}{x \, 3^n}  \int^{x}_{1}   \hspace{-0.15cm} \, \frac{dt}{  t } \sum^{n-1}_{k = 0} \frac{ 3^k  \ln^k(t)  }{k!}  +  \frac{1}{x} \int^{x}_{1}   \hspace{-0.25cm} dt \, \overline{\phi_n}(t)   \hspace{0.5cm} \\
& = & \frac{ 1 }{ x } \big(\frac{ 1 }{ 2^{n +1 } } - \frac{ 1 }{ 3^{n +1 } } \big) + \frac{1}{x \, 2^{n+1}} \sum^{n-1}_{k = 0} \frac{  2^{k+1}  \ln^{k+1}(x)  }{(k+1)!} \hspace{0.5cm}  \\
& &   - \frac{1}{x \, 3^{n+1}}   \sum^{n-1}_{k = 0} \frac{  3^{k+1}  \ln^{k+1}(x)  }{(k+1)!}  +  \overline{\phi_{n+1}}(x) \hspace{1cm}  
\end{eqnarray} \QEDA

Let's prove the $4^{th}$ point. We proceed by recurrence. For $n=1$, we retrieve the definition of $\overline{\phi_{1}}(x)$. Let's assume that it is true up to $n$ and let's prove it for $n+1$.
We have thanks to the integral order change:
\begin{eqnarray}  
\overline{\phi}_{n+1}(x) & = &  \frac{1}{x} \int^{x}_{1} dt\, \overline{\phi}_{n}(t) \\
& = &    \frac{1}{(n-1)!}   \frac{1}{x} \int^{x}_{1} \frac{dt}{t} \int^{t}_{1} ds\, \overline{\phi}_{0}(s)  \big( \ln(\frac{t}{s}) \big)^{n-1}  \hspace{0.5cm} \\
& = &   \frac{1}{(n-1)!}   \frac{1}{x} \int^{x}_{1} ds\, \overline{\phi}_{0}(s)   \int^{x}_{s}  \frac{dt}{t} \big( \ln(\frac{t}{s}) \big)^{n-1}   \hspace{0.5cm} \\
& = & \frac{1}{n!}  \frac{1}{x} \int^{x}_{1} ds\, \overline{\phi}_{0}(s) \big( \ln(\frac{x}{s}) \big)^{n} 
\end{eqnarray}
And this proves our point of the lemma. \QEDA  \\
Let's now prove the fifth point. Let's take $x \geq 1$. 
Thanks to d'Alembert's criterion, we have for each $s \in [1, x]$, $ \lim_{n \to +\infty }  \frac{ 2^n \, \big( \ln(\frac{x}{s}) \big)^{n}  }{n!}   =  0 $. From the point $(4)$ of this lemma, we apply the dominated convergence theorem to prove that:
\begin{eqnarray}  \lim_{n \to +\infty } 2^n  \,   \overline{\phi}_{n}(x) & = &   \frac{1}{x} \int^{x}_{1} \lim_{n \to +\infty }  \frac{  2^n  \big( \ln(\frac{x}{s}) \big)^{n}  }{n!}  \overline{\phi}_{0}(s) ds  = 0 \end{eqnarray} 
From point $(3)$, we can conclude that:

\begin{eqnarray} \hspace{-1cm}   
\lim_{n \to +\infty } 2^n  \,  \phi_{n}(x) & = & \lim_{n \to +\infty } \Bigg( 2^n  \,   \overline{\phi}_{n}(x)  + \frac{1}{ \, x } \Bigg[ \sum^{n-1}_{k = 0} \frac{ 2^k  \ln^k(x)  }{k!} - \big(\frac{2}{3}\big)^n \sum^{n-1}_{k = 0} \frac{ 3^k  \ln^k(x)  }{k!} \Bigg]   \Bigg)   \hspace{0.7cm}  \\
& = & \lim_{n \to +\infty }  \frac{1}{ \, x } \Bigg[ \sum^{n-1}_{k = 0} \frac{ \ln^k(x^2)  }{k!} - \big(\frac{2}{3}\big)^n \sum^{n-1}_{k = 0} \frac{ \ln^k(x^3)  }{k!} \Bigg]   \hspace{0.7cm}  \\
& = &  \frac{1}{ \, x } \exp({\ln(x^2)})  = x \hspace{0.7cm} 
 \end{eqnarray}  
For $0 \leq x < 1 $, we have from equation (63), $ \phi_{n}(x) = \frac{x}{2^n} - \frac{x^2}{3^n} $, for each $n \geq 0$. Therefore $\lim_{n \to +\infty } 2^n  \,  \phi_{n}(x)  = x$. Hence the proof of the point (5). Let's now prove the point (6). We have for each $x >0$:
\begin{eqnarray} 0 \leq \phi_0(x) \leq x \end{eqnarray} 
Therefore by integrating the equation above we get:
\begin{eqnarray}  0 \leq \frac{1}{x} \int^{x}_{0}   \hspace{-0.25cm} dt \, \phi_0(t)  \leq \frac{1}{x} \int^{x}_{0}   \hspace{-0.25cm} dt \, t   \end{eqnarray} 
Therefore
\begin{eqnarray} 0 \leq \phi_1(x) \leq \frac{ x}{2} \end{eqnarray} 
By recurrence, we easily conclude that for each $n \geq 1$: 
\begin{eqnarray} 0 \leq \phi_n(x) \leq \frac{ x}{2^n} \end{eqnarray} 
\QEDA

Let's now prove the last point. We use Taylor's Theorem with integral form of remainder applied on the exponential function $x \to e^{x}$. 

For $x \geq 0$ and $n \geq 1$:
\begin{eqnarray}  e^{x} =  \sum^{n-1}_{k = 0} \frac{ x^k  }{k!} +  \frac{ 1 }{(n-1)!}  \int^{x}_{0} dt\,  e^{t}  \big( x - t \big)^{n-1}   \end{eqnarray} 
So, let's take $x \geq 1$. We can write:
\begin{eqnarray}  e^{\ln(x^2)} =  \sum^{n-1}_{k = 0} \frac{ \big( \ln(x^2)  \big)^k  }{k!} +  \frac{ 1 }{(n-1)!}  \int^{\ln(x^2)}_{0} dt\,  e^{t}  \big( \ln(x^2) - t \big)^{n-1} 
 \end{eqnarray} 
We do a change of variable and write for $\ln(x^2)$:
\begin{eqnarray}  x^2 =  \sum^{n-1}_{k = 0} \frac{ 2^k  \ln^k(x)  }{k!} +  \frac{ 2^{n} }{(n-1)!}  \int^{x}_{1} dt\,  t  \big( \ln(\frac{x}{t}) \big)^{n-1}
 \end{eqnarray} 
Therefore
\begin{eqnarray}  \frac{x}{2^n} = \frac{1}{x \, 2^n} \sum^{n-1}_{k = 0} \frac{ 2^k  \ln^k(x)  }{k!} +  \frac{ 1 }{x \, (n-1)!}  \int^{x}_{1} dt\,  t  \big( \ln(\frac{x}{t}) \big)^{n-1} 
 \end{eqnarray} 
And the same for  $\ln(x^3)$:
\begin{eqnarray}  \frac{x^2}{3^n} = \frac{1}{x \, 3^n} \sum^{n-1}_{k = 0} \frac{ 3^k  \ln^k(x)  }{k!} +  \frac{ 1 }{x \, (n-1)!}  \int^{x}_{1} dt\,  t^2  \big( \ln(\frac{x}{t}) \big)^{n-1} 
 \end{eqnarray} 
Therefore
\begin{eqnarray} \hspace{-1cm}
 \frac{x}{2^n} - \frac{x^2}{3^n}  & = & \frac{1}{ 2^n \, x } \sum^{n-1}_{k = 0} \frac{ 2^k  \ln^k(x)  }{k!} - \frac{1}{ 3^n \, x } \sum^{n-1}_{k = 0} \frac{ 3^k  \ln^k(x)  }{k!}  +   \frac{ 1 }{x \, (n-1)!}  \int^{x}_{1} dt\,  (t - t^2) \big( \ln(\frac{x}{t}) \big)^{n-1} \hspace{0.7cm}
\end{eqnarray} \QEDA

\end{proof}

\begin{lemma}
Let's consider two variables $\sigma$ and $\tau$ such that $0 < \sigma \leq \frac{1}{2}$ and $\tau >0$ such that $s = \sigma + \textit{i} \tau$ is a zeta zero. 
Let’s define the sequence of functions $E_{n,\sigma}(x)$, $F_{n,\sigma}(x)$ and $G_{n,\sigma}(x)$ over $[1, +\infty)$ for each $n \geq 0$ as follows:

\begin{eqnarray}    \hspace{-0.7cm}  
 E_{n,\sigma}(x) & = &   \int^{x}_{1}   \hspace{-0.25cm} dt  \,  \frac{ \cos (\tau \ln{(t )})  }{  \, t^{2 + \sigma }  }  \, \phi_{n}(t)    \hspace{0.7cm}  \\
 F_{n,\sigma}(x) & = &   \int^{x}_{1}   \hspace{-0.25cm} dt  \,  \frac{ \sin (\tau \ln{(t )})  }{  \, t^{2 + \sigma }  }  \, \phi_{n}(t)    \hspace{0.7cm} \\
 G_{n,\sigma}(x) & = &   E_{n,\sigma}(x) F_{n, 1- \sigma}(x)   -   E_{n,1 - \sigma}(x) F_{n,  \sigma}(x)  \hspace{0.7cm}  
\end{eqnarray}

\begin{enumerate}

\item There exists $a_0 > 1$, there is an integer $n_0$ such that for each $n \geq n_0$:
\begin{eqnarray}    \hspace{-0.7cm}  
  G_{n,\sigma}(a_0) & > & 0   \hspace{0.7cm}     \hspace{0.7cm}  
\end{eqnarray}

\item For big enough $n$ there exists $x_n > a_0 $ such as:
\begin{eqnarray}    \hspace{-0.7cm}  
  G_{n,\sigma}(x_n) & = & 0   \hspace{0.7cm}     \hspace{0.7cm}  
\end{eqnarray}

\item If the sequence $(x_n) $ is unbounded then:
\begin{eqnarray}    \hspace{-0.7cm}  
  \sigma & = & \frac{1}{2}   \hspace{0.7cm}     \hspace{0.7cm}  
\end{eqnarray}

\item If the sequence $(x_n) $ is bounded then:
\begin{eqnarray}    \hspace{-0.7cm}  
  \sigma & = & \frac{1}{2}   \hspace{0.7cm}     \hspace{0.7cm}  
\end{eqnarray}

\end{enumerate}

\end{lemma}

\begin{proof}

Let's prove the first point. Let's prove it by contradiction. So let's assume that the opposite is true. Therefore for each $a_0 > 1$, for each integer $n_0$, there exists $n \geq n_0$ such that:
\begin{eqnarray}    \hspace{-0.7cm}  
  G_{n,\sigma}(a_0) & \leq & 0   \hspace{0.7cm}     \hspace{0.7cm}  
\end{eqnarray}
So, let's take $a_0 > 1$. Therefore for each $n \geq 1$, there exists $k_n \geq n$ such that:  
\begin{eqnarray}    \hspace{-0.7cm}  
  G_{k_n,\sigma}(a_0+\frac{1}{n}) & \leq & 0   \hspace{0.7cm}     \hspace{0.7cm}  
\end{eqnarray} 
By construction we have: \begin{eqnarray}  \lim_{n \to +\infty } k_{n} & = & +\infty \end{eqnarray} 
We apply the dominated convergence theorem on the sequences of functions: 
\begin{eqnarray} 
f^1_n(x) & = &  \frac{ \cos (\tau \ln{(t )})  }{  \, t^{2 + \sigma }  }  \, 2^n\phi_{n}(t)      \hspace{1.cm} \\  
f^2_n(x) & = &  \frac{ \sin (\tau \ln{(t )})  }{  \, t^{2 + \sigma }  }  \, 2^n\phi_{n}(t)    \hspace{1.cm} \\  
f^3_n(x) & = &  \frac{ \cos (\tau \ln{(t )})  }{  \, t^{3 - \sigma }  }  \, 2^n\phi_{n}(t)      \hspace{1.cm} \\  
f^4_n(x) & = &  \frac{ \sin (\tau \ln{(t )})  }{  \, t^{3 - \sigma }  }  \, 2^n\phi_{n}(t)    \hspace{1.cm} 
\end{eqnarray}

From lemma 1.3, we have:  \begin{eqnarray}  \lim_{n \to +\infty } 2^n  \,  \phi_{n}(x) & = & x \end{eqnarray} 
Therefore for each $x \geq 1$:

 \begin{eqnarray}  
\lim_{n \to +\infty } f^1_{k_n}(x)  & = & f_1(x)  = \frac{ \cos (\tau \ln{(x )})  }{  \, x^{ 1+\sigma }  }  \\
\lim_{n \to +\infty } f^2_{k_n}(x) & = &  f_2(x)  = \frac{ \sin (\tau \ln{(x )})  }{  \, x^{ 1+\sigma }  } \\
\lim_{n \to +\infty } f^3_{k_n}(x) & = & f_3(x)  = \frac{ \cos (\tau \ln{(x )})  }{  \, x^{ 2- \sigma }  }  \\
\lim_{n \to +\infty } f^4_{k_n}(x) & = & f_4(x)  = \frac{ \sin (\tau \ln{(x )})  }{  \, x^{ 2- \sigma }  }  
\end{eqnarray} 
Thanks to the Dominated Convergence Theorem applied on the sequence of functions in $(99-102)$ over the interval $[1,a_0]$, we apply the limit to both sides of the equation(97) as follows:
\begin{eqnarray}    \hspace{-2cm}  
\lim_{n \to +\infty }  \Big( \int^{a_0 + \frac{1}{n} }_{1}   \hspace{-0.25cm} dt  \,   f^1_n(x) \Big) \Big( \int^{a_0 + \frac{1}{n} }_{1}   \hspace{-0.25cm} dt  \,   f^4_n(x)  \Big) & \leq &  \lim_{n \to +\infty }   \Big( \int^{a_0 + \frac{1}{n} }_{1}   \hspace{-0.25cm} dt  \,   f^2_n(x)  \Big) \Big( \int^{a_0 + \frac{1}{n} }_{1}   \hspace{-0.25cm} dt  \,   f^3_n(x) \Big)  \hspace{0.8cm} 
\end{eqnarray}

\begin{eqnarray}    \hspace{-2cm}  
\Big( \int^{a_0}_{1}   \hspace{-0.25cm} dt  \,  \frac{ \cos (\tau \ln{(t )})  }{  \, t^{1 + \sigma }  }  \Big) \Big( \int^{a_0}_{1}   \hspace{-0.25cm} dt  \,  \frac{ \sin (\tau \ln{(t )})  }{  \, t^{2 - \sigma }  }  \Big) & \leq &   \Big( \int^{a_0}_{1}   \hspace{-0.25cm} dt  \,  \frac{ \cos (\tau \ln{(t )})  }{  \, t^{2 - \sigma }  }  \Big) \Big( \int^{a_0}_{1}   \hspace{-0.25cm} dt  \,  \frac{ \sin (\tau \ln{(t )})  }{  \, t^{1 + \sigma }  }   \Big)  \hspace{0.8cm} 
\end{eqnarray}
From the lemma 1.1 result, we conclude:
\begin{eqnarray}    \hspace{-0.7cm}  
g(a_0) & = & J(a_0, 1+\sigma, \tau) I(a_0, 2-\sigma, \tau)  - J(a_0, 2-\sigma, \tau) I(a_0, 1+\sigma, \tau) \leq 0  \hspace{1.2cm}  
\end{eqnarray}
Where the function $g$ is defined as follows: \begin{eqnarray}    \hspace{-2.5cm}  
g(x) & = & J(x, 1+\sigma, \tau) I(x, 2-\sigma, \tau)  - J(x, 2-\sigma, \tau) I(x, 1+\sigma, \tau)   \hspace{0.8cm}  \\
& = & K(1+\sigma, \tau) K(2-\sigma, \tau)  \Bigg( \Big(1 + \frac{\sigma(1-\sigma)}{\tau^2}\Big) \frac{ \sin (\tau \ln{(x )})  }{  \, x^{\sigma }  } \Big( 1 - \frac{1}{x^{1-2\sigma}} \Big)  \hspace{0.8cm}   \\
&   & + \frac{(1-2\sigma)}{\tau}\Bigg[  \frac{ \cos (\tau \ln{(x )}) }{ x^{\sigma } }  \Big(1 + \frac{1}{x^{1-2\sigma}} \Big) -  \Big( 1 + \frac{1}{x} \Big) \Bigg] \Bigg) \hspace{0.8cm} 
\end{eqnarray}
Therefore for each $x \geq 1$: \begin{eqnarray}    \hspace{-0.0cm}  
g(x) & \leq & 0  \hspace{1.2cm}  
\end{eqnarray}
Which is a contradiction since the function $g$ oscillates between negative and positive values over $[1,+\infty)$. In fact, from lemma 2.3 below, there exists a constant $\tau_0 \in [1.8549, 1.8554]$ such that for each $\tau \geq \tau_0$, we have $g(e^{\frac{5 \pi}{2\tau}}) \geq 0$ for each $\sigma \leq \frac{1}{2}$. And for  $\tau < \tau_0$, there is no zeta zero.

Therefore there is $a_0$ such that $g(a_0) > 0$ and hence the first point of the lemma is proved. \QEDA   \\
Let's now prove the $2^{nd}$ point of the lemma. From the lemma 1.3, we can write:
\begin{eqnarray}    \hspace{-1cm}  
  \lim_{x \to +\infty } G_{n,\sigma}(x) & = &   \lim_{x \to +\infty } E_{n,\sigma}(x) F_{n, 1- \sigma}(x)   -   E_{n,1 - \sigma}(x) F_{n,  \sigma}(x)  \hspace{0.7cm}  \\
& = & \Big( \int^{+\infty}_{1}   \hspace{-0.25cm} dt  \,  \frac{ \cos (\tau \ln{(t )})  }{  \, t^{2 + \sigma }  }  \, \phi_{n}(t) \Big) \Big( \int^{+\infty}_{1}   \hspace{-0.25cm} dt  \,  \frac{ \sin (\tau \ln{(t )})  }{  \, t^{3 - \sigma }  }  \, \phi_{n}(t) \Big)  \\
&  & -   \Big( \int^{+\infty}_{1}   \hspace{-0.25cm} dt  \,  \frac{ \cos (\tau \ln{(t )})  }{  \, t^{3 - \sigma }  }  \, \phi_{n}(t) \Big) \Big( \int^{+\infty}_{1}   \hspace{-0.25cm} dt  \,  \frac{ \sin (\tau \ln{(t )})  }{  \, t^{2 + \sigma }  }  \, \phi_{n}(t) \Big)  \hspace{0.7cm} 
\end{eqnarray}
Therefore, we can write:
\begin{eqnarray}    \hspace{-1cm}  
  \lim_{x \to +\infty } G_{n,\sigma}(x)  & = & - \frac{\tau(1-2\sigma)}{ \norm{s(1-s)}^2 } \big(  \frac{1}{ 2^n} - \frac{1}{ 3^n} \big)^2  < 0 \hspace{2cm}  
\end{eqnarray}

Since the function $G_{n,\sigma}$ is continuous over $[1, +\infty)$. From (92) and $(118)$ and thanks to the Mean value theorem, we can conclude that there exists an $x_n > a_0$ such that:
\begin{eqnarray}    \hspace{0cm}  
   G_{n,\sigma}(x_n) & = & 0 \hspace{0.7cm} 
\end{eqnarray} \QEDA

Let's prove the $3^{rd}$ point. So, let's assume the sequence $(x_n) $ is unbounded. There exists a subsequence\footnote{Similar to Bolzano–Weierstrass theorem in the case of a bounded sequence.} $(x_{\lambda(n)})$ that tends to infinity $+\infty$. Without loss of generality, we assume that the limit of the sequence $(x_n)$ is $+\infty$.
From the equation (93), we can deduce that:
\begin{eqnarray}    \hspace{-2cm}  
\Big( \int^{x_n}_{1}   \hspace{-0.25cm} dt  \,  \frac{ \cos (\tau \ln{(t )})  }{  \, t^{2 + \sigma }  }  \, \phi_{n}(t) \Big) \Big( \int^{x_n}_{1}   \hspace{-0.25cm} dt  \,  \frac{ \sin (\tau \ln{(t )})  }{  \, t^{3 - \sigma }  }  \, \phi_{n}(t) \Big) & = &   \Big( \int^{x_n}_{1}   \hspace{-0.25cm} dt  \,  \frac{ \cos (\tau \ln{(t )})  }{  \, t^{3 - \sigma }  }  \, \phi_{n}(t) \Big) \Big( \int^{x_n}_{1}   \hspace{-0.25cm} dt  \,  \frac{ \sin (\tau \ln{(t )})  }{  \, t^{2 + \sigma }  }  \, \phi_{n}(t) \Big)  \hspace{0.8cm} 
\end{eqnarray}
We multiply the two side of the equation above by $4^n$, we can write:
\begin{eqnarray}    \hspace{-2.5cm}  
\Big( \int^{x_n}_{1}   \hspace{-0.25cm} dt  \,  \frac{ \cos (\tau \ln{(t )})  }{  \, t^{2 + \sigma }  }  \, 2^n\phi_{n}(t) \Big) \Big( \int^{x_n}_{1}   \hspace{-0.25cm} dt  \,  \frac{ \sin (\tau \ln{(t )})  }{  \, t^{3 - \sigma }  }  \, 2^n\phi_{n}(t) \Big) & = &   \Big( \int^{x_n}_{1}   \hspace{-0.25cm} dt  \,  \frac{ \cos (\tau \ln{(t )})  }{  \, t^{3 - \sigma }  }  \, 2^n\phi_{n}(t) \Big) \Big( \int^{x_n}_{1}   \hspace{-0.25cm} dt  \,  \frac{ \sin (\tau \ln{(t )})  }{  \, t^{2 + \sigma }  }  \, 2^n\phi_{n}(t) \Big)  \hspace{0.8cm} 
\end{eqnarray}
We apply the Dominated Convergence Theorem on the sequences of functions $f^1_n, f^2_n, f^3_n, f^4_n$ defined in the equations $(99-102)$ since we have the dominance condition for the functions $f^i_{1 \leq i \leq 4}$: 
\begin{eqnarray} 
\left| f^i_n(x)_{1 \leq i \leq 2} \right| & \leq &  \left|  \frac{ 1  }{  \, t^{2 + \sigma }  } \, 2^n\phi_{n}(t)   \right| \leq    \frac{ 1 }{  \, t^{1 + \sigma }  }  \hspace{1.cm}  \\
\left| f^i_n(x)_{3 \leq i \leq 4} \right| & \leq &  \left|  \frac{ 1  }{  \, t^{3 - \sigma }  } \, 2^n\phi_{n}(t)   \right| \leq    \frac{ 1 }{  \, t^{2 - \sigma }  }  \hspace{1.cm}  
\end{eqnarray}

Thanks to the Dominated Convergence Theorem applied over the interval $[1,+\infty)$, we apply the limit to both sides of the equation(121) as follows:
\begin{eqnarray}    \hspace{-2cm}  
\lim_{n \to +\infty }  \Big( \int^{x_n}_{1}   \hspace{-0.25cm} dt  \,   f^1_n(x) \Big) \Big( \int^{x_n}_{1}   \hspace{-0.25cm} dt  \,   f^4_n(x)  \Big) & = &  \lim_{n \to +\infty }   \Big( \int^{x_n}_{1}   \hspace{-0.25cm} dt  \,   f^2_n(x)  \Big) \Big( \int^{x_n}_{1}   \hspace{-0.25cm} dt  \,   f^3_n(x) \Big)  \hspace{0.7cm} 
\end{eqnarray}
Therefore
\begin{eqnarray}    \hspace{-2cm}  
  \Big( \int^{+\infty}_{1}   \hspace{-0.25cm} dt  \, \lim_{n \to +\infty }   f^1_n(x) \Big) \Big( \int^{+\infty}_{1}   \hspace{-0.25cm} dt  \, \lim_{n \to +\infty }   f^4_n(x)  \Big) & = &  \lim_{n \to +\infty }   \Big( \int^{+\infty}_{1}   \hspace{-0.25cm} dt  \,   \lim_{n \to +\infty } f^2_n(x)  \Big) \Big( \int^{+\infty}_{1}   \hspace{-0.25cm} dt  \,  \lim_{n \to +\infty }  f^3_n(x) \Big)  \hspace{0.8cm} 
\end{eqnarray}
Therefore
\begin{eqnarray}    \hspace{-2cm}  
\Big( \int^{+\infty}_{1}   \hspace{-0.25cm} dt  \,  \frac{ \cos (\tau \ln{(t )})  }{  \, t^{1 + \sigma }  }  \Big) \Big( \int^{+\infty}_{1}   \hspace{-0.25cm} dt  \,  \frac{ \sin (\tau \ln{(t )})  }{  \, t^{2 - \sigma }  }  \Big) & = &   \Big( \int^{+\infty}_{1}   \hspace{-0.25cm} dt  \,  \frac{ \cos (\tau \ln{(t )})  }{  \, t^{2 - \sigma }  }  \Big) \Big( \int^{+\infty}_{1}   \hspace{-0.25cm} dt  \,  \frac{ \sin (\tau \ln{(t )})  }{  \, t^{1 + \sigma }  }   \Big)  \hspace{0.8cm} 
\end{eqnarray}
From the lemma 1.1 result, we conclude:
\begin{eqnarray}    \hspace{-0cm}  
\frac{\sigma}{\tau} K(1+\sigma, \tau) K(2-\sigma, \tau)  & = &   \frac{1 - \sigma}{\tau} K(2-\sigma, \tau) K(1+\sigma, \tau)  \hspace{0.7cm} 
\end{eqnarray}
Therefore: \begin{eqnarray}    \hspace{-0cm}  
\sigma  & = &   1 - \sigma = \frac{1}{2} \hspace{0.7cm} 
\end{eqnarray} \QEDA \\
Let's now prove the last point of the lemma. So, let's assume the sequence $(x_n) $ is bounded by an upper bound $A > 0$. Thanks to Bolzano–Weierstrass theorem, there exists a subsequence $(x_{\lambda(n)})$ that converges to a finite limit $a \geq a_0 > 1$. And without loss of generality, we can assume that the limit of the sequence $(x_n)$ is $a$.

Let's define the functions $g_1$, $g_2$ and $g_3$ as follows:
\begin{eqnarray}    \hspace{-0cm}  
g_1(x) & = & J(x, 1+\sigma, \tau) I(x, 2-\sigma, \tau)  - J(x, 2-\sigma, \tau) I(x, 1+\sigma, \tau)   \hspace{0.7cm}  \\ 
g_2(x) & = & J(x, 1+\sigma, \tau) I(x, 1-\sigma, \tau)  + J(x, \sigma, \tau) I(x, 2-\sigma, \tau)  \\
&  & - J(x, 2-\sigma, \tau) I(x, \sigma, \tau) -  J(x, 1-\sigma, \tau) I(x, 1+\sigma, \tau)  \hspace{0.7cm}  \\
g_3(x) & = & J(x, \sigma, \tau) I(x, 1-\sigma, \tau)  - J(x, 1-\sigma, \tau) I(x, \sigma, \tau)   \hspace{0.7cm}  
\end{eqnarray}

From the lemma 1.3 point 3) and point 7), we can write for each $x \geq 1$:
\begin{eqnarray}  \phi_{n}(x) =  \psi_{n}(x)  + \underbrace{ \overline{\phi}_{n}(x) - \overline{\psi}_{n}(x) }_{  \delta_n(x)   }  \end{eqnarray} 
Where
\begin{eqnarray}  
\psi_{n}(x)  & = & \frac{x}{2^n} -  \frac{x^2}{3^n} \\
 \overline{\phi_{n}}(x) - \overline{\psi}_{n}(x)  & = &  \frac{1}{(n-1)!}  \frac{1}{x}    \int^{x}_{1} dt\,  \big(\overline{\phi_{0}}(t) + t^2 - t \big)  \big( \ln(\frac{x}{t}) \big)^{n-1}   \hspace{1cm}
\end{eqnarray}
Therefore for each $x \in [1, A]$:
\begin{eqnarray}  
 \left| \delta_n(x)  = \overline{\phi_{n}}(x) - \overline{\psi}_{n}(x)  \right| & \leq &  \frac{A^2}{(n-1)!}  \frac{\big( \ln(x) \big)^{n-1} (x-1) }{x}  \leq    \frac{ A^2  \big( \ln(A) \big)^{n-1} }{(n-1)!}        \hspace{1cm}
\end{eqnarray}

We inject the equations $(133-134)$ into the equation $(120)$ to write the following. For big enough $n$:
\begin{eqnarray}    \hspace{-3cm}  
 \Big( \int^{x_n}_{1}   \hspace{-0.25cm} dt  \,  (\psi_{n}(x)  + \delta_n(x) \frac{ f_1(x)}{x} \Big) \Big( \int^{x_n}_{1}   \hspace{-0.25cm} dt  \,   (\psi_{n}(x)  +\delta_n(x) \frac{ f_4(x)}{x}   \Big)  =   \Big( \int^{x_n}_{1}   \hspace{-0.25cm} dt  \,  (\psi_{n}(x)  + \delta_n(x) \frac{ f_2(x)}{x}  \Big) \Big( \int^{x_n}_{1}   \hspace{-0.25cm} dt  \,   (\psi_{n}(x)  + \delta_n(x) \frac{ f_3(x)}{x}  \Big)  \hspace{0.1cm} 
\end{eqnarray}
Therefore:
\begin{eqnarray}    \hspace{-3cm}  
  \int^{x_n}_{1}   \hspace{-0.25cm} dt  \,  \psi_{n}(x)  \frac{ f_1(x)}{x} \int^{x_n}_{1}   \hspace{-0.25cm} dt  \,   \psi_{n}(x) \frac{ f_4(x)}{x}  -   \int^{x_n}_{1}   \hspace{-0.25cm} dt  \,  \psi_{n}(x)  \frac{ f_2(x)}{x}  \int^{x_n}_{1}   \hspace{-0.25cm} dt  \,   \psi_{n}(x)  \frac{ f_3(x)}{x}   \hspace{0.3cm} \\
=    \int^{x_n}_{1}   \hspace{-0.25cm} dt  \,  \delta_n(x)  \frac{ f_2(x)}{x}  \int^{x_n}_{1}   \hspace{-0.25cm} dt  \,   \phi_{n}(x) \frac{ f_3(x)}{x} +  \int^{x_n}_{1}   \hspace{-0.25cm} dt  \,  \delta_n(x) \frac{ f_3(x)}{x}  \int^{x_n}_{1}   \hspace{-0.25cm} dt  \,   \phi_{n}(x) \frac{ f_2(x)}{x}  \hspace{0.3cm} \\
- \int^{x_n}_{1}   \hspace{-0.25cm} dt  \,  \delta_n(x)  \frac{ f_1(x)}{x}  \int^{x_n}_{1}   \hspace{-0.25cm} dt  \,   \phi_{n}(x) \frac{ f_4(x)}{x} +  \int^{x_n}_{1}   \hspace{-0.25cm} dt  \,  \delta_n(x)  \frac{ f_4(x)}{x}  \int^{x_n}_{1}   \hspace{-0.25cm} dt  \,   \phi_{n}(x) \frac{ f_1(x)}{x}  \hspace{0.3cm}
\end{eqnarray}
Thanks to (136), from the equations (138-140) we can write:
\begin{eqnarray}    \hspace{-1cm}  
  \int^{x_n}_{1}   \hspace{-0.25cm} dt  \,  \psi_{n}(x)  \frac{ f_1(x)}{x} \int^{x_n}_{1}   \hspace{-0.25cm} dt  \,   \psi_{n}(x) \frac{ f_4(x)}{x}  -   \int^{x_n}_{1}   \hspace{-0.25cm} dt  \,  \psi_{n}(x)  \frac{ f_2(x)}{x}  \int^{x_n}_{1}   \hspace{-0.25cm} dt  \,   \psi_{n}(x)  \frac{ f_3(x)}{x}  =    \mathcal{O}( \frac{ \big( \ln(A) \big)^{n-1} }{(n-1)!}   )  \hspace{0.3cm}
\end{eqnarray}
We inject the equation (134) into the equation $(141)$. After simplification, we write the following:
\begin{eqnarray}    \hspace{-1cm}  
 \frac{1}{4^n} g_1(x_n) - \frac{1}{6^n} g_2(x_n) + \frac{1}{9^n} g_3(x_n)  =    \mathcal{O}( \frac{ \big( \ln(A) \big)^{n-1} }{(n-1)!}   )  \hspace{0.3cm}
\end{eqnarray}
Therefore
\begin{eqnarray}    \hspace{-1cm}  
 \big(\frac{3}{2}\big)^{2n} g_1(x_n) -  \big(\frac{3}{2}\big)^{n} g_2(x_n)  +  g_3(x_n)  =    \mathcal{O}( \frac{ \big( 9 \ln(A) \big)^{n-1} }{(n-1)!}   )  \hspace{0.3cm}
\end{eqnarray}

Since the limit of the sequence $( \frac{ \big( 9 \ln(A) \big)^{n-1} }{(n-1)!}   )$ is zero thanks to Stirling's formula. Since the functions $g_1$, $g_2$ and $g_3$ are continuous and therefore bounded over the interval $[1,A]$. There exists $n_0$ such that for each $n \geq n_0$, we have:
\begin{eqnarray}    \hspace{-1cm}  
g_1(x_n)  & = & 0    \hspace{0.3cm} \\
g_2(x_n)  & = & 0   \hspace{0.3cm} 
\end{eqnarray}
Therefore
\begin{eqnarray}  
g_1(a)   =  g_2(a)  =  g_3(a)   = 0 \hspace{0.3cm} 
\end{eqnarray}

Hence $\sigma$ must be equal to $1-\sigma$. 
Therefore $\sigma = \frac{1}{2}$. \QEDA \\ 
This ends the proof of the Riemann Hypothesis. 

\end{proof}


\subsection{Conclusion}


We saw that if $s$ is a zeta zero, then real part $\Re(s)$ can only be $\frac{1}{2}$. Therefore the Riemann's Hypothesis is true: \textit{The non-trivial zeros of $\zeta(s)$ have real part equal to $\frac{1}{2}$}. In the next article, we will apply the same method to prove the Generalized Riemann Hypothesis (GRH). 

\subsection*{Acknowledgments}
I would like to thank Farhat Latrach, Giampiero Esposito, Jacques Gélinas, Michael Milgram, Léo Agélas, Ronald F. Fox, Kim Y.G, Masumi Nakajima, Maksym Radziwill and Shekhar Suman for thoughtful comments and discussions on my paper versions on the RH. All errors are mine.

\section{Appendix}

\begin{lemma}
Let's consider two variables $\sigma$ and $\tau$ such that $\sigma >0$ and $\tau >0$. Let's define two integrals $I(+, \sigma, \tau)$ and $I(-, \sigma, \tau)$ as follows: 
\begin{eqnarray} I(+, \sigma, \tau) & = &  \int^{+\infty}_{1}  \frac{  \cos^{+} (\tau \ln(x)) }{ x^{1 + \sigma} }   dx  \\
I(-, \sigma, \tau) & = &  \int^{+\infty}_{1}  \frac{  \cos^{-} (\tau \ln(x)) }{ x^{1 + \sigma} }   dx  
\end{eqnarray}  
Therefore
\begin{eqnarray} I(+, \sigma, \tau) & = &  K(\sigma, \tau)   \frac{\sigma}{\tau} +    K(\sigma, \tau)  \frac{  e^{-\frac{\pi \sigma}{2\tau}}   } { 1 - e^{-\frac{\pi \sigma}{\tau}}  }      \hspace{1cm} \\
I(-, \sigma, \tau) & = &       K(\sigma, \tau)     \frac{  e^{-\frac{\pi \sigma}{2\tau}}  }{ 1 - e^{-\frac{\pi \sigma}{\tau}}  }    \hspace{1cm} 
\end{eqnarray}
Where \begin{eqnarray} K(\sigma, \tau) & = &   \frac{\tau}{\sigma^2 + \tau^2}   \hspace{1cm} \\
 \cos^{+} (x )  & = & \max \big( \cos (x ), 0  \big) \hspace{1cm}\\
 \cos^{-} (x)  & = & \max \big( - \cos(x), 0  \big)   \hspace{1cm}
 \end{eqnarray} 
\end{lemma}

\begin{proof}

So, let's start.

\begin{eqnarray} \hspace{-2cm}  I(+, \sigma, \tau)  & = &    \int^{+\infty}_{1}   \hspace{-0.25cm} dx \frac{  \cos^{+} (\tau \ln{(x )})  }{  \, x^{1 + \sigma }  }  \\
& = & \int^{  e^{\frac{\pi}{2\tau}}  }_{1}   \hspace{-0.25cm} dx \frac{  \cos (\tau \ln{(x )})  }{  \, x^{1 + \sigma }  }  + \sum^{+\infty}_{k = 0}   \int^{e^{\frac{(4k+5)\pi}{2\tau}}}_{e^{\frac{(4k+3)\pi}{2\tau}}}   \hspace{-0.25cm} dx \frac{   \cos (\tau \ln{x})  }{  \, x^{1 + \sigma }  }    \hspace{1cm} \\
& = & \int^{  e^{\frac{\pi}{2\tau}}  }_{1}   \hspace{-0.25cm} dx \frac{  \cos (\tau \ln{(x )})  }{  \, x^{1 + \sigma }  }  +  \sum^{+\infty}_{k = 0}  I(k, \sigma, \tau)  \hspace{1cm}  \end{eqnarray} Where 
\begin{eqnarray} \hspace{-2cm}   I(+, \sigma, \tau)  & = &   \sum^{+\infty}_{k = 0}   \int^{e^{\frac{(4k+5)\pi}{2\tau}}}_{e^{\frac{(4k+3)\pi}{2\tau}}}   \hspace{-0.25cm} dx \frac{   \cos (\tau \ln{x})  }{  \, x^{1 + \sigma }  }      \hspace{1cm}    \end{eqnarray} 

From lemma 1.1 we can write the following:
\begin{eqnarray} \hspace{-2cm}   I(+, \sigma, \tau)  & = &   \sum^{+\infty}_{k = 0}   \int^{e^{\frac{(4k+5)\pi}{2\tau}}}_{e^{\frac{(4k+3)\pi}{2\tau}}}   \hspace{-0.25cm} dx \frac{   \cos (\tau \ln{x})  }{  \, x^{1 + \sigma }  }      \hspace{1cm} \\
& = & K(\sigma, \tau)  \Bigg(  e^{-\frac{(4k+5)\pi \sigma}{2\tau}} + e^{-\frac{(4k+3)\pi \sigma}{2\tau}}     \Bigg)     \end{eqnarray} 

Therefore
\begin{eqnarray} \hspace{-2cm}  I(+, \sigma, \tau)  & = &   K(\sigma, \tau)   \Bigg(  \frac{\sigma}{\tau} +   e^{-\frac{\pi \sigma}{2\tau}} \Bigg)   +  K(\sigma, \tau)   \sum^{+\infty}_{k = 0}  \Bigg(  e^{-\frac{(4k+5)\pi \sigma}{2\tau}} + e^{-\frac{(4k+3)\pi \sigma}{2\tau}}     \Bigg)   \hspace{1cm}  \\
& = &  K(\sigma, \tau)   \Bigg(  \frac{\sigma}{\tau} +  e^{-\frac{\pi \sigma}{2\tau}} \Bigg)   +  K(\sigma, \tau)   \sum^{+\infty}_{k = 0}  \Bigg( e^{-\frac{5\pi \sigma}{2\tau}}  e^{-\frac{2k\pi \sigma}{\tau}} + e^{-\frac{3\pi \sigma}{2\tau}} e^{-\frac{2k\pi \sigma}{\tau}}     \Bigg)   \hspace{1cm}  \\
& = &  K(\sigma, \tau)   \Bigg(  \frac{\sigma}{\tau} +   e^{-\frac{\pi \sigma}{2\tau}} \Bigg)   +  K(\sigma, \tau)   \frac{ e^{-\frac{5\pi \sigma}{2\tau}}  + e^{-\frac{3\pi \sigma}{2\tau}} } { 1 - e^{-\frac{2\pi \sigma}{\tau}}  }     \hspace{1cm}  \\
& = &    K(\sigma, \tau)   \Bigg(  \frac{\sigma}{\tau} +   e^{-\frac{\pi \sigma}{2\tau}} \Bigg)   +  K(\sigma, \tau)   e^{-\frac{3\pi \sigma}{2\tau}}  \frac{ 1  + e^{-\frac{\pi \sigma}{\tau}} } { 1 - e^{-\frac{2\pi \sigma}{\tau}}  }     \hspace{1cm} \\
& = &    K(\sigma, \tau)   \Bigg(  \frac{\sigma}{\tau} +   e^{-\frac{\pi \sigma}{2\tau}} \Bigg)   +  K(\sigma, \tau)   e^{-\frac{3\pi \sigma}{2\tau}}  \frac{ 1 } { 1 - e^{-\frac{\pi \sigma}{\tau}}  }     \hspace{1cm} \\ 
& = & K(\sigma, \tau)   \frac{\sigma}{\tau} +    K(\sigma, \tau)   e^{-\frac{\pi \sigma}{2\tau}}    \Bigg( 1 +   \frac{  e^{-\frac{\pi \sigma}{\tau}}  } { 1 - e^{-\frac{\pi \sigma}{\tau}}  } \Bigg)     \hspace{1cm}  \\
& =  & K(\sigma, \tau)   \frac{\sigma}{\tau} +  K(\sigma, \tau)   \frac{  e^{-\frac{\pi \sigma}{2\tau}}   } { 1 - e^{-\frac{\pi \sigma}{\tau}}  }    \end{eqnarray}

And therefore

\begin{eqnarray} \hspace{-2cm}  I(-, \sigma, \tau)  & = &  I(+\infty, \sigma, \tau) - I(+, \sigma, \tau)  \\  
& = &   K(\sigma, \tau)    \frac{  e^{-\frac{\pi \sigma}{2\tau}}   } { 1 - e^{-\frac{\pi \sigma}{\tau}}  }   \hspace{1cm} \end{eqnarray}

\QEDA

\end{proof}

\begin{lemma}
Let's consider $f$ a continuous function over $[1, +\infty)$. Let's $\phi$ be a non-null positive function such that $f \phi$ and $\phi$ are integrable functions over $[1, +\infty)$ with:
\begin{eqnarray}  0 < \int^{+\infty}_{1}   \hspace{-0.25cm} dx \, \phi(x)  & < &  +\infty   \end{eqnarray}
And
 \begin{eqnarray}   \int^{+\infty}_{1}   \hspace{-0.25cm} dx \, \phi(x) \, f(x)  & < & +\infty   \end{eqnarray} 
Therefore, there exists a $c \in (1, +\infty)$ such that: 
\begin{eqnarray}   \int^{+\infty}_{1}   \hspace{-0.25cm} dx \, \phi(x) \, f(x)  & = & f(c)   \int^{+\infty}_{1}   \hspace{-0.25cm} dx \, \phi(x)   \end{eqnarray} 
\end{lemma}
\begin{proof}
Let's define the real $\lambda$ as following:
\begin{eqnarray}   \lambda  & = &  \frac{ \int^{+\infty}_{1}   dx \, \phi(x) \, f(x)  } {   \int^{+\infty}_{1}  dx \, \phi(x) }   \end{eqnarray} 
We have by construction that:
\begin{eqnarray}   \int^{+\infty}_{1}   dx \, \phi(x) \, \big( f(x) - \lambda \big)  & = & 0  \end{eqnarray} 
Therefore, if for each $x > 1$, we have $ f(x) > \lambda$, then,we will have:
\begin{eqnarray}   \int^{+\infty}_{1}   dx \, \phi(x) \, \big( f(x) - \lambda \big)  & > & 0  \end{eqnarray}
Which is a contradiction. We will reach a similar contradiction if we assume $ f(x) < \lambda$ for each $x > 1$. Therefore, there exists $c \in (1, +\infty)$ such that $f(c) = \lambda$. \QEDA
\end{proof}

\begin{lemma}
Let's consider two variables $\sigma$ and $\tau$ such that $0 < \sigma < 1$ and $\tau >0$ and $s = \sigma + \textit{i} \tau$ is a zeta zero.
Therefore:  
 \begin{eqnarray} \tau  &  >  & \frac{\sqrt{3 \pi}}{2} \sqrt{ 1 + \frac{1}{1 + \frac{1}{2} \sqrt{\frac{\pi}{3}} } } \sim 1.9786    \end{eqnarray} 
\end{lemma}

\begin{proof}
From the lemma 1.3, we have:
\begin{eqnarray}  \hspace{-1cm}  \int^{+\infty}_{1}   \hspace{-0.25cm} dx \frac{ \phi(x)  \cos (\tau \ln{(x )})  }{  \, x^{2 + \sigma }  }   & = &    \frac{\sigma }{  \sigma^2  + \tau^2   }   +  \frac{ \big(1 - \sigma \big) }{  (1-\sigma)^2  + \tau^2  } \hspace{1cm}  \\
 \int^{+\infty}_{1}   \hspace{-0.25cm} dx \frac{ \phi(x)  \sin (\tau \ln{(x )})  }{  \, x^{2 + \sigma }  }   & = &    \frac{\tau }{  \sigma^2  + \tau^2   }   -  \frac{ \tau }{  (1-\sigma)^2  + \tau^2  } \hspace{1cm}  
 \end{eqnarray}  

Therefore, we can write the following:
\begin{eqnarray}  \hspace{-0cm}  \int^{+\infty}_{1}   \hspace{-0.25cm} dx \frac{ \phi(x)  \cos (\tau \ln{(x )})    }{  \, x^{2 + \sigma }  }   & = &       \int^{+\infty}_{1}   \hspace{-0.25cm} dx \frac{ \phi(x)  \cos^{+} (\tau \ln{(x )})  }{  \, x^{2 + \sigma }  }  -  \int^{+\infty}_{1}   \hspace{-0.25cm} dx \frac{ \phi(x)  \cos^{-} (\tau \ln{(x )})  }{  \, x^{2 + \sigma }  } \hspace{1cm}
 \end{eqnarray}  
Where \begin{eqnarray}    \cos^{+} (x )  & = & \max \big( \cos (x ), 0  \big) \\
 \cos^{-} (x)  & = & \max \big( - \cos(x), 0  \big)   \hspace{1cm}
 \end{eqnarray}  

From lemma 2.2, there exists $c_1 > 1 $ and $c_2 > 1$ such that:
\begin{eqnarray}  \hspace{-0cm}  \int^{+\infty}_{1}   \hspace{-0.25cm} dx \frac{ \phi(x)  \cos (\tau \ln{(x )})  }{  \, x^{2 + \sigma }  }   & = &    \frac{ \phi(c_1) }{c^{2\sigma}_{1}}   \int^{+\infty}_{1}   \hspace{-0.25cm} dx \frac{   \cos^{+} (\tau \ln{(x )})  }{  \, x^{2 - \sigma }  }  - \frac{\phi(c_2)}{c^{2\sigma}_{2}}  \int^{+\infty}_{1}   \hspace{-0.25cm} dx \frac{  \cos^{-} (\tau \ln{(x )})  }{  \, x^{2 - \sigma }  } \hspace{1cm}
 \end{eqnarray}  
Let's denote $\alpha_1 = \frac{ \phi(c_1) }{c^{2\sigma}_{1}}  $ and $\alpha_2 =  \frac{\phi(c_2)}{c^{2\sigma}_{2}} $. 
We have $0 < \alpha_1 < \frac{1}{4}$ and $ 0 < \alpha_2 < \frac{1}{4}$.

\paragraph{Case 1: $\alpha_1 <= \alpha_2$ } In this case we can write $ \alpha_1 = \alpha_2 - \epsilon $ with $0 < \epsilon < \frac{1}{4}$.
Therefore we can write from the equation (176) that:
\begin{eqnarray}  \hspace{-1cm}     \alpha_1 \int^{+\infty}_{1}   \hspace{-0.25cm} dx \frac{   \cos^{+} (\tau \ln{(x )})  }{  \, x^{2 - \sigma }  }  - \alpha_2  \int^{+\infty}_{1}   \hspace{-0.25cm} dx \frac{  \cos^{-} (\tau \ln{(x )})  }{  \, x^{2 - \sigma }  } & = &    \frac{\sigma }{  \sigma^2  + \tau^2   }   +  \frac{ \big(1 - \sigma \big) }{  (1-\sigma)^2  + \tau^2  } \hspace{1cm} \hspace{1cm}
 \end{eqnarray}  

Therefore
\begin{eqnarray}  \hspace{-1cm}     \alpha_1 \int^{+\infty}_{1}   \hspace{-0.25cm} dx \frac{   \cos^{+} (\tau \ln{(x )})  }{  \, x^{2 - \sigma }  }  - \alpha_1  \int^{+\infty}_{1}   \hspace{-0.25cm} dx \frac{  \cos^{-} (\tau \ln{(x )})  }{  \, x^{2 - \sigma }  } -  \epsilon  \int^{+\infty}_{1}   \hspace{-0.25cm} dx \frac{  \cos^{-} (\tau \ln{(x )})  }{  \, x^{2 - \sigma }  } & = &    \frac{\sigma }{  \sigma^2  + \tau^2   }   +  \frac{ \big(1 - \sigma \big) }{  (1-\sigma)^2  + \tau^2  } \hspace{1cm}
 \end{eqnarray}  
Therefore

\begin{eqnarray}  \hspace{-1cm}     \alpha_1 \int^{+\infty}_{1}   \hspace{-0.25cm} dx \frac{   \cos (\tau \ln{(x )})  }{  \, x^{2 - \sigma }  }  -  \epsilon  \int^{+\infty}_{1}   \hspace{-0.25cm} dx \frac{  \cos^{-} (\tau \ln{(x )})  }{  \, x^{2 - \sigma }  } & = &    \frac{\sigma }{  \sigma^2  + \tau^2   }   +  \frac{ \big(1 - \sigma \big) }{  (1-\sigma)^2  + \tau^2  } \hspace{1cm}
 \end{eqnarray}  

Therefore

\begin{eqnarray}  \hspace{-1cm}     \alpha_1 \frac{ \big(1 - \sigma \big) }{  (1-\sigma)^2  + \tau^2  }  -  \epsilon  \int^{+\infty}_{1}   \hspace{-0.25cm} dx \frac{  \cos^{-} (\tau \ln{(x )})  }{  \, x^{2 - \sigma }  } & = &    \frac{\sigma }{  \sigma^2  + \tau^2   }   +  \frac{ \big(1 - \sigma \big) }{  (1-\sigma)^2  + \tau^2  } \hspace{1cm} \end{eqnarray} 
Therefore
\begin{eqnarray} 
 -  \epsilon  \int^{+\infty}_{1}   \hspace{-0.25cm} dx \frac{  \cos^{-} (\tau \ln{(x )})  }{  \, x^{2 - \sigma }  } & = &    \frac{\sigma }{  \sigma^2  + \tau^2   }   + (1 - \alpha_1)  \frac{ \big(1 - \sigma \big) }{  (1-\sigma)^2  + \tau^2  } \hspace{1cm}
 \end{eqnarray}  
Which is a contradiction.

\paragraph{Case 2: $\alpha_1 > \alpha_2$ } In this case we can write $ \alpha_1 = \alpha_2 + \epsilon $ with $0 < \epsilon < \frac{1}{4}$.
Therefore we can write from the equation (176) that:
\begin{eqnarray}  \hspace{-1cm}     \alpha_1 \int^{+\infty}_{1}   \hspace{-0.25cm} dx \frac{   \cos^{+} (\tau \ln{(x )})  }{  \, x^{2 - \sigma }  }  - \alpha_2  \int^{+\infty}_{1}   \hspace{-0.25cm} dx \frac{  \cos^{-} (\tau \ln{(x )})  }{  \, x^{2 - \sigma }  } & = &    \frac{\sigma }{  \sigma^2  + \tau^2   }   +  \frac{ \big(1 - \sigma \big) }{  (1-\sigma)^2  + \tau^2  } \hspace{1cm} \hspace{1cm}
 \end{eqnarray}  

Therefore
\begin{eqnarray}  \hspace{-1cm}     \alpha_2 \int^{+\infty}_{1}   \hspace{-0.25cm} dx \frac{   \cos^{+} (\tau \ln{(x )})  }{  \, x^{2 - \sigma }  }  - \alpha_2  \int^{+\infty}_{1}   \hspace{-0.25cm} dx \frac{  \cos^{-} (\tau \ln{(x )})  }{  \, x^{2 - \sigma }  } +  \epsilon  \int^{+\infty}_{1}   \hspace{-0.25cm} dx \frac{  \cos^{+} (\tau \ln{(x )})  }{  \, x^{2 - \sigma }  } & = &    \frac{\sigma }{  \sigma^2  + \tau^2   }   +  \frac{ \big(1 - \sigma \big) }{  (1-\sigma)^2  + \tau^2  } \hspace{1cm}
 \end{eqnarray}  
Therefore

\begin{eqnarray}  \hspace{-1cm}     \alpha_2 \int^{+\infty}_{1}   \hspace{-0.25cm} dx \frac{   \cos (\tau \ln{(x )})  }{  \, x^{2 - \sigma }  }  +  \epsilon  \int^{+\infty}_{1}   \hspace{-0.25cm} dx \frac{  \cos^{+} (\tau \ln{(x )})  }{  \, x^{2 - \sigma }  } & = &    \frac{\sigma }{  \sigma^2  + \tau^2   }   +  \frac{ \big(1 - \sigma \big) }{  (1-\sigma)^2  + \tau^2  } \hspace{1cm}
 \end{eqnarray}  

Therefore

\begin{eqnarray}  \hspace{-1cm}     \alpha_2 \frac{ \big(1 - \sigma \big) }{  (1-\sigma)^2  + \tau^2  }  +  \epsilon  \int^{+\infty}_{1}   \hspace{-0.25cm} dx \frac{  \cos^{+} (\tau \ln{(x )})  }{  \, x^{2 - \sigma }  } & = &    \frac{\sigma }{  \sigma^2  + \tau^2   }   +  \frac{ \big(1 - \sigma \big) }{  (1-\sigma)^2  + \tau^2  } \hspace{1cm} \end{eqnarray} 
Therefore
\begin{eqnarray} 
 + \epsilon  \int^{+\infty}_{1}   \hspace{-0.25cm} dx \frac{  \cos^{+} (\tau \ln{(x )})  }{  \, x^{2 - \sigma }  } & = &    \frac{\sigma }{  \sigma^2  + \tau^2   }   + (1 - \alpha_2)  \frac{ \big(1 - \sigma \big) }{  (1-\sigma)^2  + \tau^2  } \hspace{1cm}
 \end{eqnarray}  

From lemma 2.1, we have:
\begin{eqnarray} 
 \int^{+\infty}_{1}   \hspace{-0.25cm} dx \frac{  \cos^{+} (\tau \ln{(x )})  }{  \, x^{2 - \sigma }  } & = &  \frac{ \big(1 - \sigma \big) }{  (1-\sigma)^2  + \tau^2  } +  \frac{ \tau }{  (1-\sigma)^2  + \tau^2  }   \frac{  e^{-\frac{\pi (1-\sigma)}{2\tau} }   }{ 1 - e^{-\frac{\pi (1-\sigma)}{\tau} } }   \hspace{1cm}
 \end{eqnarray}  

Therefore
\begin{eqnarray} 
   \frac{  \epsilon \, \tau }{  (1-\sigma)^2  + \tau^2  }   \frac{  e^{-\frac{\pi (1-\sigma)}{2\tau} }   }{ 1 - e^{-\frac{\pi (1-\sigma)}{\tau} } }     & = &    \frac{\sigma }{  \sigma^2  + \tau^2   }   + (1 - \alpha_1)  \frac{ \big(1 - \sigma \big) }{  (1-\sigma)^2  + \tau^2  } \hspace{1cm}
 \end{eqnarray}  
Therefore
\begin{eqnarray} 
   \epsilon  \, \tau  \frac{  e^{-\frac{\pi (1-\sigma)}{2\tau} }   }{ 1- e^{-\frac{\pi (1-\sigma)}{\tau} } }    & = & \sigma    \frac{  (1-\sigma)^2  + \tau^2  }{  \sigma^2  + \tau^2   }   + (1 - \alpha_1)  \big(1 - \sigma \big) \hspace{1cm}
 \end{eqnarray} 

Therefore
\begin{eqnarray} 
   \epsilon  \, \tau  \frac{  e^{-\frac{\pi (1-\sigma)}{2\tau} }   }{ 1- e^{-\frac{\pi (1-\sigma)}{\tau} } }    & = &    \frac{  \sigma (1-2\sigma)  }{  \sigma^2  + \tau^2   }   + \sigma +  \big(1 - \sigma \big)  - \alpha_1  \big(1 - \sigma \big)\hspace{1cm}
 \end{eqnarray} 

We have $0 < \sigma < \frac{1}{2}$. Therefore

\begin{eqnarray} 
   \epsilon  \, \tau  \frac{  e^{-\frac{\pi (1-\sigma)}{2\tau} }   }{ 1 - e^{-\frac{\pi (1-\sigma)}{\tau} } }     & > & 1 -  \alpha_1  \big(1 - \sigma \big) \hspace{1cm}
 \end{eqnarray} 

We have $0 < \alpha_1 < \frac{1}{4}$ and $1 > 1 - \sigma > \frac{1}{2}$. Therefore:  $1 -  \alpha_1  \big(1 - \sigma \big)  > \frac{3}{4}$. 
Therefore
\begin{eqnarray} 
   \epsilon  \, \tau   \frac{  e^{-\frac{\pi (1-\sigma)}{2\tau} }   }{ 1 - e^{-\frac{\pi (1-\sigma)}{\tau} } }     & > & \frac{3}{4}  \hspace{1cm}
 \end{eqnarray} 

But we also have $0 < \epsilon < \frac{1}{4}$. Therefore

\begin{eqnarray} 
    \frac{3}{4}   & < &  \frac{\tau}{4}  \frac{  e^{-\frac{\pi (1-\sigma)}{2\tau} }   }{ 1 - e^{-\frac{\pi (1-\sigma)}{\tau} } }     \hspace{1cm}
 \end{eqnarray} 

Since $e^x \geq 1 + x$ for each $x \geq 0$, therefore:

\begin{eqnarray} 
    \frac{\tau}{3}   & > &      e^{\frac{\pi (1-\sigma)}{2\tau} }   -  e^{-\frac{\pi (1-\sigma)}{2\tau} }    \\
& > & 1 + \frac{\pi (1-\sigma)}{2\tau} -  \frac{1}{  1 + \frac{\pi (1-\sigma)}{2\tau}   }   \hspace{1cm}
 \end{eqnarray} 
Since the function $x \to \frac{x}{2\tau + x}$ is increasing over $[0,+\infty)$, we can write:
\begin{eqnarray} 
    \frac{4 \tau^2}{3 \pi}   & > &   1 + \frac{4 \tau}{4\tau + \pi}     \hspace{1cm}
 \end{eqnarray} 
Since $\tau >0$, therefore, 
\begin{eqnarray} 
    \tau & > & \frac{ \sqrt{3 \pi} }{2}      \hspace{1cm}
 \end{eqnarray} 

Since the function $x \to \frac{x}{x +\pi}$ is increasing over $[0,+\infty)$. Therefore
\begin{eqnarray} 
    \frac{4 \tau^2}{3 \pi}   & > &   1 + \frac{  2\sqrt{3 \pi} }{ 2 \sqrt{3 \pi}   + \pi}     \hspace{1cm}
 \end{eqnarray} 

Therefore
\begin{eqnarray} \tau  &  >  & \frac{\sqrt{3 \pi}}{2} \sqrt{ 1 + \frac{1}{1 + \frac{1}{2} \sqrt{\frac{\pi}{3}} } }    \end{eqnarray} 
\remark{We can reiterate the same procedure using the equation (201) to improve the minimum bound further. The limit bound $\tau_{min}$ is actually the root of the equation $\frac{16}{3\pi}\tau^3 + \frac{4}{3 } \tau^2 - 8 \tau - \pi = 0$. $\tau_{min} \sim 2.01271781$. } \QEDA 

\end{proof}  

\begin{lemma}
Let's consider two variables $\sigma$ and $\tau$ such that $0 < \sigma < 1$ and $\tau >0$. Let's consider the function $g$ defined in the equation $(111)$.
Therefore there exists a constant $\tau_0 \in [1.8549, 1.8554]$ such that:
\begin{enumerate}
\item For $\tau \geq \tau_0$, we have $g(e^{\frac{5 \pi}{2\tau}}) \geq 0$ for each $\sigma \leq \frac{1}{2}$. 
\item For  $\tau < \tau_0$, $s = \sigma + \textit{i} \tau$ cannot be a zeta zero.
\end{enumerate}
\end{lemma}

\begin{proof}

Let's prove the first point of the lemma. We will ignore $K(1+\sigma, \tau) K(2-\sigma, \tau) $ from the function $g$ expression as follows:
\begin{eqnarray}    \hspace{-0cm}  
g(x) & = &  \Big(1 + \frac{\sigma(1-\sigma)}{\tau^2}\Big) \frac{ \sin (\tau \ln{(x )})  }{  \, x^{\sigma }  } \Big( 1 - \frac{1}{x^{1-2\sigma}} \Big)  \hspace{0.8cm}   \\
&   & + \frac{(1-2\sigma)}{\tau}\Bigg[  \frac{ \cos (\tau \ln{(x )}) }{ x^{\sigma } }  \Big(1 + \frac{1}{x^{1-2\sigma}} \Big) -  \Big( 1 + \frac{1}{x} \Big) \Bigg]  \hspace{0.8cm} 
\end{eqnarray}
 We can find $a$ such that $g(x) > 0$ of the form $a = e^{\frac{\lambda \pi}{\tau}}$ where $\lambda$ is of the form $\lambda(k) = 2k + \frac{1}{2}$ and $k \in \mathbb{N}$. For such $x$ we have:
\begin{eqnarray}    \hspace{-0.0cm}  
g(a)  & = & \Big(1 + \frac{\sigma(1-\sigma)}{\tau^2}\Big) e^{-\frac{\lambda \pi \sigma}{\tau}} \Big( 1 - e^{-\frac{\lambda \pi (1-2\sigma)}{\tau}} \Big)  - \frac{(1-2\sigma)}{\tau} \Big( 1 + e^{-\frac{\lambda \pi}{\tau}} \Big)  \hspace{0.8cm} 
\end{eqnarray}

Let's denote $x = \frac{\lambda \pi}{\tau}$ and $f_{\lambda}(x) = g(e^x)$. Therefore
\begin{eqnarray}    \hspace{-0.0cm}  
 f_{\lambda}(x,\sigma)  & = & \Big(1 + \frac{\sigma(1-\sigma)}{ \big(\lambda \pi \big)^2} x^2 \Big) \Big( e^{- x \sigma } - e^{-{ x (1-\sigma)}}  \Big)  - \frac{(1-2\sigma)}{\lambda \pi} x \Big( 1 + e^{-x} \Big)  \hspace{0.8cm} 
\end{eqnarray}

We can use the asymptotic expansion around 0 ($x \to 0$) of $f_{\lambda}$ to get an idea of what is happening around 0. For $\epsilon >0$ small enough, we can write: \begin{eqnarray}    \hspace{0.0cm}  
f_{\lambda}(1+\epsilon, \sigma) & = & \big(1 - 2\sigma \big) \big(1- \frac{2}{\lambda \pi} \big) \Big( \epsilon  - \frac{1}{2}  \epsilon^2 \Big) + \mathcal{O}(\epsilon^3) \hspace{0.3cm}  
\end{eqnarray}

To get $f_{\lambda}(1+\epsilon, \sigma) > 0$, we need $\lambda \geq \frac{2}{\pi}$ and therefore $k \geq 1$ since $1- 2\sigma >0$. For the rest of this section we will take $\lambda = \lambda(1) = \frac{5}{2}$. We will plot the function $f_{\frac{5}{2}}$ over $[0,5]$ for different values of $\sigma$ from $0$ to $\frac{1}{2}$.
\begin{figure}
  \includegraphics[width=\linewidth]{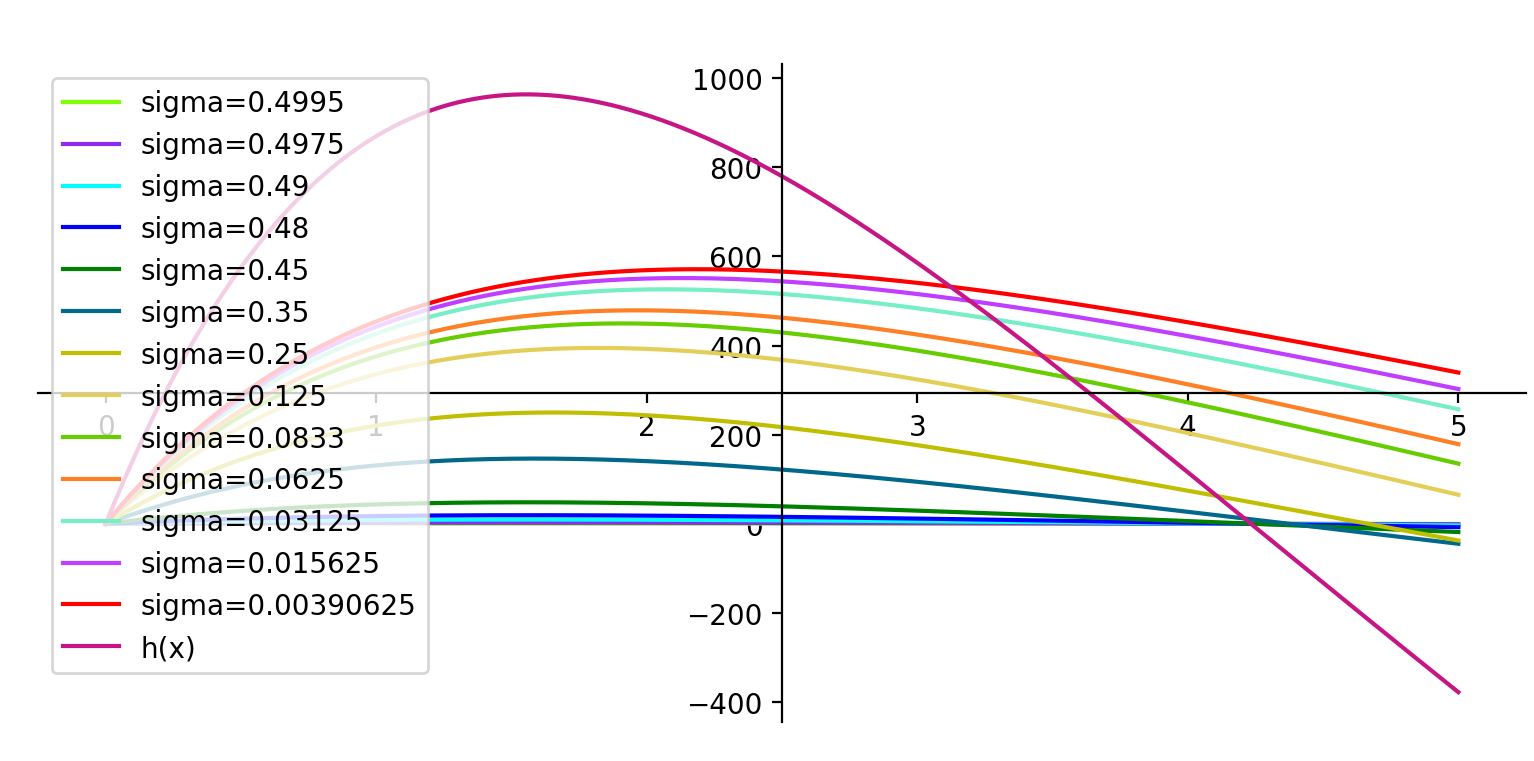}
  \caption{Functions $f_{\frac{5}{2}}(x, \sigma)$ for different values of $\sigma$ from $0$ to $\frac{1}{2}$ and function $h_{\frac{5}{2}}(x)$ over $[0,5]$ .   }
  \label{fig:boat1}
\end{figure}

From figure $1$, there exists a constant $x_0$ such that for each $x \in [0, x_0]$, we have $f_{\frac{5}{2}}(x) \geq 0$ for each $\sigma \leq \frac{1}{2}$. Also, the more $\sigma$ is close to $\frac{1}{2}$, the more $f_{\lambda}(x_0, \sigma)$ is close to zero.  The asymptotic expansion of $f_{\lambda}$ when $\sigma$ tends to $\frac{1}{2}$ is as follows: \begin{eqnarray}    \hspace{0.0cm}  
f_{\lambda}(x, \sigma) & = & \big(1 - 2\sigma \big) \underbrace{ x e^{-\frac{x}{2}} \Big( 1 + \frac{ x^2 }{(2\lambda \pi)^2} - \frac{ \big( e^{\frac{x}{2}} + e^{-\frac{x}{2}} \big) }{\lambda \pi} \Big) }_{ h_{\lambda}(x) } + \mathcal{O}(( \sigma - \frac{1}{2} )^3) \hspace{0.8cm}  
\end{eqnarray}
$x_0$ is actually the first zero of the function $h_{\lambda}(x)$ that is different from zero.  Numerical tests showed that $x_0 \in (4.233, 4.234)$. Therefore,  there exists a constant $\tau_0 = \frac{5 \pi}{2 x_0}$ with $\tau_0 \in [1.8549, 1.8554]$ such that for each $\tau \geq \tau_0$, we have $g(e^{\frac{5 \pi}{2\tau}}) \geq 0$ for each $\sigma \leq \frac{1}{2}$. \QEDA \\

Let's prove the second point of the lemma.  Thanks to lemma 2.3, we have $\tau > \tau_1 = \frac{\sqrt{3 \pi}}{2} \sqrt{ 1 + \frac{1}{1 + \frac{1}{2} \sqrt{\frac{\pi}{3}} } } \sim 1.9786 $. Therefore if $\tau < \tau_0$, then $\tau < \tau_1$, therefore $s = \sigma + \it{i} \tau$ cannot be a zeta zero.  \QEDA

\end{proof}


\begin{thebibliography}{9}

\bibitem{RiemannMain} 
Bernhard Riemann.
\textit{On the Number of Prime Numbers less than a Given Quantity}
\\\texttt{\url{https://www.claymath.org/sites/default/files/ezeta.pdf}}

\bibitem{ZetaFullBook} 
Aleksandar Ivic. 
\textit{The Riemann Zeta-Function: Theory and Applications}

\bibitem{ZetaFullBook} 
Peter Borwein, Stephen Choi, Brendan Rooney, and Andrea Weirathmueller
\textit{The Riemann Hypothesis: A Resource for the Afficionado and Virtuoso Alike}
\\\texttt{\url{http://wayback.cecm.sfu.ca/~pborwein/TEMP\_PROTECTED/book.pdf}}

 
\bibitem{VeisdalBlog} 
J\o rgen Veisdal.
\textit{The Riemann Hypothesis, explained}
\\\texttt{\url{https://medium.com/cantors-paradise/the-riemann-hypothesis-explained-fa01c1f75d3f}}


\bibitem{ThaiPham} 
Thai Pham.
\textit{Dirichlet’s Theorem on Arithmetic Progressions}
\\\texttt{https://web.stanford.edu/~thaipham/papers/MIT\_18.104\_Review\_Paper.pdf}


\bibitem{Hardy} 
 G. H. Hardy.
\textit{The general theory of dirichlet series.}
\\\texttt{https://archive.org/details/generaltheoryofd029816mbp/page/n9}

\bibitem{Garrett} 
 Garrett, Paul.
\textit{Primes in arithmetic progressions, 2011.}
\\\texttt{http : //www.math.umn.edu/~garrett/m/mfms/notes\_c/dirichlet.pdf}



\bibitem{Habil} 
 Eissa D. Habil.
\textit{Double Sequences and Double Series.}
\\\texttt{https://journals.iugaza.edu.ps/index.php/IUGNS/article/download/1594/1525}











\end{thebibliography}
\end{document}